\documentclass[11pt]{elsart}
\usepackage{amssymb}
\usepackage{graphicx} 
\newtheorem{theo}{Theorem}%[section]
\begin{document}

\begin{frontmatter}
\title{iso-spectral Euler-Bernoulli beams  \`a la Sophus Lie}

\author{C\'elestin Wafo Soh\thanksref{w}}
%\date{}
%\maketitle
%\begin{center}
\address{Mathematics Department,
 College of Science, Engineering, and Technology  Jackson State
 University,
 JSU Box 17610, 1400 J R Lynch St., \\Jackson, MS 39217, USA}
\ead{wafosoh@yahoo.com, celestin.wafo@jsums.edu}
 \thanks[w]{To the loving  memory of my brother L\'eopold Fotso Simo.}

\begin{abstract}
We obtain iso-spectral Euler-Bernoulli beams by using factorization
 and Lie symmetry techniques. The canonical Euler-Bernoulli
beam operator is factorized as the product of a second-order linear
differential operator and its adjoint. The factors are then reversed
to obtain iso-spectral beams. The factorization is possible provided
the coefficients of the factors satisfy a system of non-linear
ordinary differential equations . The uncoupling of this system
yields a single non-linear third-order ordinary differential
equation. This ordinary differential equation, refer to as the {\it
principal equation}, is analyzed and solved using Lie group methods.
We show that the principal equation may admit a one-dimensional or
three-dimensional symmetry Lie algebra. When the principal system
admits a unique symmetry, the best we can do is to depress its order
by one. We obtain a one-parameter family of solutions in this case.
The maximally symmetric case is shown to be isomorphic to a Chazy
equation which is solved in closed form to derive the general
solution of the principal equation.

\end{abstract}

%\begin{keyword}
%Iso-spectral Euler-Bernoulli beams, Factorization, Lie symmetry,
%Chazy equation, Gauss hypergeometric function. \PACS
%\end{keyword}

\end{frontmatter}

\section{Introduction}

A recent study \cite{roba} suggests that efforts to model the
transverse motion of vibrating beams date back to Leonardo da Vinci.
In his discussion of the bending of beam/spring with rectangular
cross-section he wrote \cite{reti}: ``Of bending of the springs: if
a straight spring is bent, it is necessary that its convex part
become thinner and its concave part, thicker. This modification is
pyramidal, and consequently, there will never be a change in the
middle of the spring''. About hundred years after Da Vinci's attempt
to develop a beam theory, Galileo suggested an erroneous calculation
of the load carrying capacity of a transversely loaded beam. Da
Vinci  did not benefit from  Hooke's law and calculus which postdate
him whereas Galileo made the incorrect assumption that under
transverse loading, the stress is uniformly distributed on
cross-sections. The first correct and  systematic formulation of an
elasticity  theory was done by Jacob Bernoulli (1700-1782). He found
that the curvature of an  elastic beam at a given  point is
proportional to the bending moment at that same point. Taking
advantage of his uncle seminal work, Daniel Bernoulli(1707-1783)
derived the partial  differential equation that governs the motion
of a vibrating beam. Leonard Euler(1707-1705) pursued the study of
the Benoullis by using and extending their theory in his
investigation of the shape of elastic beams subjected to various
external forces. For a recent review of beam theories, the reader is
referred to the paper by Han, Benaroya and Wei  \cite{han} and
references therein.

In  Euler-Bernoulli beam theory, the transverse motion of a thin
non-uniform beam is governed by the partial differential equation
\begin{equation}
 {\partial^2 \over \partial x^2} \left [ E(x)\;I(x){\partial^2 u\over \partial x^2}\right]
 + \rho(x)\; \alpha(x){\partial^2 u\over \partial t^2}=0,\quad 0\le x\le L, \label{i1}
\end{equation}
where $E$ is the modulus of elasticity (a.k.a Young's modulus),
$\alpha(x)$ is  the cross-sectional area at $x$, $I(x)$ is the area
moment of inertia of the cross-section located at $x$ about a
neutral axis, $\rho$ is the density, and $L$ is the length of the
beam. Equation (\ref{i1}) is solved subject to appropriate boundary
conditions at $x=0$ and $x=L$. Some common boundary conditions and
their physical meaning are
\begin{eqnarray}
& & {\partial^2 u \over \partial x^2}=0,\;u=0 \mbox{ ({\small hinged
end}) }, \label{i2} \\
& &      {\partial u \over \partial
x}=0,\; u=0 \mbox{ ({\small clamped end})},  \label{i2'}  \\
 & & {\partial^2 u \over \partial x^2}=0,\;
 {\partial\over \partial x}\left ( E(x)I(x){\partial^2 u \over
\partial x^2}\right )=0
 \mbox{ ({\small free end}) }, \label{i2''}\\
 && {\partial u \over \partial
x}=0,\;{\partial\over \partial x}\left ( E(x)I(x){\partial^2 u \over
\partial x^2}\right )=0 \mbox{ ({\small sliding end})}\label{i2'''}  .
\end{eqnarray}
 The ansatz $u(t,x)=\e^{i\omega t}Y(x)$, where $\omega$ is the
frequency of vibrations, solves Eq. (\ref{i1}) provided
\begin{equation}
{d^2\over dx^2}\left ( f(x){d^2Y\over dx^2}\right
)=\omega^2m(x)Y,\quad 0\le x\le L, \label{i3}
\end{equation}
where
\begin{equation}
 f(x)=E(x)\;I(x),\quad m(x)=\rho(x)\;\alpha(x). \label{i4}
\end{equation}

Under Barcilon's  transformation \cite{barc,gott}
\begin{equation}
z=\int_0^x (m/f)^{1/4},\quad Y=(m^3f)^{-1/8}\; U, \label{i5}
\end{equation}
Equation (\ref{i5}) becomes
\begin{equation}
{d^4U\over dz^2}+{d\over dz}\left [ A(z){dU \over dz}\right
]+B(z)U=\omega^2U, \quad 0\le z\le  l, \label{i6}
\end{equation}
where the correct expression of $A$ and $B$ are given by Gottlieb
\cite{gott}. Save for the clamped end conditions,  Barcilon's
transformation leaves the other boundary  conditions invariant under
some constraints given explicitly  in \cite{gott}. As we shall see
in the sequel, the canonical form (\ref{i6}) plays a crucial in
inverse spectral problems associated with Euler-Bernoulli beams.

In many practical problems such as the structural identification
through non-destructive means, one is interested in determining the
physical characteristics of the beam i.e. $f(x)$ and $m(x)$ from few
frequencies (spectra). That is, one tries to `hear' the shape of the
beam. Since the vibration of uniform beams (i.e. $f$ and  $m$ are
constant functions) is well understood, an interesting problem
related to the inverse spectral problem is to find non-unform beams
that have the same length and share the same spectra with the
uniform beam under the same boundary conditions. This problem was
comprehensively studied by Gottlieb \cite{gott} using the canonical
form (\ref{i6}). He found seven classes of clamped  beams that are
iso-spectral with the standard unit-coefficient beam ($A=0=B$). Some
of these beams where rediscovered by Abrate \cite{abra} (see
\cite{gotl} for comments) using a different paradigm. For discrete
beams, Gladwell \cite{glad}, introduced two procedures for finding
iso-spectral families. His first procedure relies on shifted
orthogonal matrix triangulation  whereas the second utilizes Toda
flow.

This paper is dedicated to the search of iso-spectral beams using
Lie symmetry methods together with a  factorization method initiated
by P\"ochel and Trubowitz \cite{potr} for Sturm-Liouville operator.
To be more precise, the factorization method for the canonical beam
operator leads to a system of nonlinear ordinary differential
equation that we  shall solve by Lie's reduction method. We have
organized the present paper as follows. There are four section of
which this introduction is the first. In section 2, we briefly
motivate the use of factorization in the search of iso-spectral
beams. The factorization approach leads to a system nonlinear
ordinary differential equations. We employ Lie's symmetry theory to
study its continuous symmetries in Section 3. Section 4 employs the
symmetries calculated in Section 3 to find exact solutions of the
system of nonlinear ordinary differential equations derived in
Section 2. This results in new families of iso-spectral beams.

\section{Iso-spectral deformation through factorization}
The main idea is to factorize the canonical beam operator as the
product of a linear second-order differential operator and its
adjoint, then swap the order of the operators to obtain iso-spectral
operators. This idea was already used for vibrating rods by P\"ochel
and Trubowitz \cite{potr}. Ghanbari \cite{ghan} attempted the
factorization method on vibrating beams but failed to solve the
resulting system of nonlinear ordinary differential equations.

Let us define the beam operator as the differential operator
\begin{equation}
 \mathfrak{L}={d^4\over dz^4}+{d\over dz} \left (A(z){d\over dz}
 \right )+ B(z). \label{f1}
\end{equation}
 Factorize  the differential operator $\mathfrak{L}$ as
\begin{equation}
\mathfrak{L}=\mathfrak{R}^*\mathfrak{R}, \label{f2}
\end{equation}
where
\begin{equation}
\mathfrak{R}={d^2\over dz^2}+r(z){d\over dz}+ s(z), \quad
\mathfrak{R}^*={d^2\over dz^2}-{d\over dz} [r(z)\;.\;]+ s(z),
\label{f3}
\end{equation}
$r$ and $s$ are smooth functions of their argument, and
$\mathfrak{R}^*$ is the adjoint of $\mathfrak{R}$. Simple reckoning
shows that the factorization given in Eq. (\ref{f2}) constraints the
functions $r$ and $s$ to satisfy the following system of nonlinear
ordinary differential equations \cite{ghan}.
\begin{eqnarray}
r'-r^2 +2s & =& A(z) \label{f3} \\
s'' -(rs)'+s^2 &=& B(z), \label{f4}
\end{eqnarray}
where the prime stands for differentiation with respect to $z$.

Reversing the operators in the factorization  (\ref{f2}), we
obtained the operator \cite{ghan}
\begin{equation}
\hat{\mathfrak{L}}= \mathfrak{R}\mathfrak{R}^* = {d^4\over
dz^4}+{d\over dz} \left (\hat A(z){d\over dz}
 \right )+ \hat B(z), \label{f5}
\end{equation}
where
\begin{eqnarray}
\hat A(z) & = & 2s-3r'-r^2, \label{f6}\\
\hat B(z) &=& s^2+s''-r'''-rr''+rs'-sr'. \label{f7}
\end{eqnarray}

Define an {\it eigenpair} of $\mathfrak{L}$ as being a pair
$(\omega,U)$, where $U$ is not the zero function,  such that
$\mathfrak{L} U =\omega^2 U$. It can be easily verified that if
$(\omega, U)$ (resp. $(\hat \omega,\hat U)$) is an eigenpair of
$\mathfrak{L}$ (resp. $\hat \mathfrak{L}$), then $(\omega,
\mathfrak{R}U)$ (resp. $(\hat\omega,\mathfrak{R}^* \hat U )$ is an
eigenpair of $\hat \mathfrak{L}$ (resp. $\mathfrak{L}$). Thus,
loosely (we have not considered the boundary conditions),
$\mathfrak{L}$ and $\hat \mathfrak{L}$ are iso-spectral. This
motivates the use of factorization in the search for iso-spectral
beam.

In order to characterize iso-spectral beams obtained via
factorization, one needs to solve Eqs (\ref{f3})-(\ref{f4}) for $r$
and $s$. Solving  Eq. (\ref{f3}) for $s$, and substituting into Eq.
(\ref{f4}) yield the single nonlinear ordinary differential equation
for $r$
\begin{equation}
r^{(3)}-3rr''-{7\over 2}r'^2+2(2r^2+A)r'-A''-r^2A-{A^2\over
2}-{r^4\over 2}+rA'+2B=0. \label{f8}
\end{equation}
In the next two sections, we shall use Lie symmetry analysis to
solve Eq. (\ref{f8}). We shall refer to Eq. (\ref{f8}) as the {\it
principal equation}.

\section{Lie symmetry analysis of the principal equation}

\label{symanal}
 A symmetry of a differential equation is an
invertible transformation of the dependents and dependents
variables that leaves the equation unchanged. The main attractive
property of symmetries is that they transform solutions into
solutions, and better, they provide a systematic route to
integration.

Historically, the Norwegian mathematician Sophus Lie (1842-1899)
realized that determining all the symmetries of a given equation
is a formidable task. However, if we restrict ourself to
symmetries that depends continuously on a small parameter and that
forms a group (called an {\it infinitesimal group} in Lie's
terminology or Lie group in modern terminology) , it is possible
to use the machinery of calculus to determine such symmetries in
an algorithmic fashion. Strikingly, Sophus Lie discovered that
symmetries are road maps to integration.

Here, we  are concerned with the determination of the Lie symmetries
of Eq. (\ref{f8}). For more details about Lie's symmetry theory, the
reader is referred to the books \cite{ovsi,blum,olve}.

A vector field
\begin{equation}
X=\xi(z,r){\partial \over \partial z}+\eta(z,r){\partial \over
\partial r}, \label{s1}
\end{equation}
 is a point symmetry of Eq. (\ref{f8}) if
  \begin{eqnarray}
     & & X^{[3]} \left[  r^{(3)}-3rr''-{7\over 2}r'^2+2(2r^2+A)r'-A''-r^2A
      \right.  \nonumber  \\
     & &\quad \quad \quad \left. \left. -{A^2\over 2}-{r^4\over 2}+rA'+2B  \right]\right|_{\mbox{Eq.}(\ref{f8})} =0.
\label{s2}
 \end{eqnarray}

 In Eq. (\ref{s2}), the notation $|_{\mbox{Eq.}(\ref{f8})} $  means that after
 expanding the left hand side, Eq. (\ref{f8}) has to be used to
 eliminate $r^{(3)}$. The operator $X^{[3]}$ is the third
 prolongation  of $X$ defined recursively by
 \begin{eqnarray}
 X^{[k+1]} &= & X^{[k]}+\eta^{[k]}{\partial \over \partial
 r^{(k)}},\quad X^{[0]}=X, \label{s3}\\
 \eta^{[k]}& = & D\left (\eta^{[k-1]}\right )+r^{(k)}D(\xi),\quad \eta^{[0]}=\eta,
 \label{s4} \\
 D &=& ={\partial \over \partial z} + r' {\partial \over \partial
 r}+ r'' {\partial \over \partial
 r'}+ \cdots . \label{s5}
 \end{eqnarray}
  Since the  symmetry coefficients $\xi$ and $\eta$ are independent $r'$ and $r''$,
  Eq. (\ref{s2}) is polynomial in these  derivatives. Thus we may  set coefficients of the monomials
  $( r')^m(r'')^n$ to zero in Eq. (\ref{s2}). It results an
  over-determined system of linear partial differential equations for $\xi$ and $\eta$ that simplifies to the following equations.

  \begin{eqnarray}
  & & \xi =  a(z),\quad \eta =-a'(z)r-2a''(z), \label{s6}\\
  & & a A' +2a' A +5a^{(3)}=0, \label{s7}\\
    & & B'+{4a'\over a}\,B={a'\over a}\,A^2+{AA'\over 2}+{a''\over
  a'}\,A'+{2a'\over a}\,A''+{2a^{(3)}\over a}\,A +{A^{(3)}\over
  2}+{a^{(5)}\over a}, \label{s8}\\
  & & 3 a' A'+2a''A+aA'' +5a^{(4)}=0, \label{s9}
  \end{eqnarray}
 where $a$ is an arbitrary smooth function of $z$. Simple calculations
 show that $${d\over dz} (\mbox{ Eq}.(\ref{s7}))\equiv \mbox{
 Eq}.(\ref{s9}).$$ Thus Eq. (\ref{s9}) is a mere differential
 consequence of Eq.(\ref{s7}). In order to solve the remaining equations, we consider the following cases.

 {\bf \tt Case I: $A=0$ and $B=0$}. In this case, we find that
 $a=k_1 z^2+k_2z+k_3$, where $k_1$, $k_2$, and $k_3$ are arbitrary
 constants. The symmetry Lie algebra is spanned by the operators
 \begin{equation}
 X_1 ={\partial \over \partial z},\quad X_2 = z{\partial\over \partial
 z}-r {\partial\over \partial r},\quad X_3 =z^2{\partial \over
 \partial z}-2(rz+2){\partial\over \partial r}\;\cdot \label{s10}
 \end{equation}
 The Lie bracket of these operators are:
 $$[X_1,X_2]=X_1,\quad [X_1,X_3]=2X_2,\quad [X_2,X_3]=X_3.$$
 The symmetry Lie
 algebra is non-solvable since its derived algebra of any order is
 nontrivial. To be more precise, if $\mathfrak{S}=<X_1,X_2,X_3>$,
 then $\mathfrak{S}^{(k)}=\mathfrak{S}$. In fact, $\mathfrak{S}\cong
 sl(2,\mathbb{R})$ \cite{pate}.

 {\bf \tt Case II: $A\ne 0$ or $B\ne0$}. Solving Eqs.
 (\ref{s7})-(\ref{s8}), we obtain after some calculations
 \begin{eqnarray}
 A & =&{5a'^2-10aa''+2C_1\over 2a^2}\;,     \label{s11}\\
 B &=& {81a'^4+12a'^2(3C_1-17aa'')+72a^2a'a^{(3)}\over 16 a^4}\nonumber  \\
 & & +{4\left (C_1^2+4C_2-6C_1aa''+21a^2a''^2-6a^2a^{(4)}\right )\over 16 a^4 }\;,     \label{s12}
 \end{eqnarray}
 where $C_1$ and $C_2$ are arbitrary integration constants. The symmetry Lie algebra is spanned by the single operator
 \begin{equation}
 \Gamma = a\;{\partial \over \partial z} - ( a'r+2a'')\;{\partial \over
 \partial r}\;\cdot \label{s13}
 \end{equation}
\section{Solutions of the principal equation and iso-spectral beams}
Here we exploit the symmetry structure of the principal equation
(\ref{f8}) to find its solutions.
\subsection{ $A$ and $B$ are given by Eqs. (\ref{s11})-(\ref{s12}) }
In this case the symmetry Lie algebra is one-dimensional. The best
we can do using Lie's integration technique is to reduce the order
of the equation by one. This is accomplished by introducing new
dependent and independent variables as independent solutions (known
as {\it basis of first order differential invariants} ) of
\begin{equation}
\Gamma^{[1]} I =0. \label{so1}
\end{equation}
A fundamental set of solutions of Eq. (\ref{so1}) are obtained by
solving the system of ordinary differential equations
\begin{equation}
{dz\over a(z)} ={dr\over -a'r-2a''}={dr'\over -2a'r'-a''r-2a^{(3)}}
\; \cdot \label{so3}
\end{equation}
The solutions of Eq. (\ref{so3})  are
\begin{equation}
a r+2a' =q_1,\quad a^2 r'+aa'r+2aa''=q_2, \label{s04}
\end{equation}
where $q_1$ and $q_2$ are arbitrary constants of integration.  Thus
we introduce the new variables
\begin{equation}
u = a r+2a',\quad v= a^2 r'+aa'r+2aa''. \label{s05}
\end{equation}
In terms of the new variables, Eq. (\ref{f8}) reads
\begin{equation}
2v^2\;{d^2v\over du^2}=6uv{dv\over du}-2v\left ({dv\over du}\right
)^2 + u^4-8u^2v+7v^2 +2C_1u^2-4C_1v-4C_2\,. \label{s06}
\end{equation}
Equation (\ref{s06}) is devoid of Lie symmetries so that we may not
further integrate Eq. (\ref{s06}) by Lie's method. Assume from now
on that $C1=0=C_2$ , and  look for a solution of Eq. (\ref{s06}) in
the form
\begin{equation}
v =ku^2, \label{s07}
\end{equation}
where $k$ is a constant to be determined. Substituting the ansatz
(\ref{s07}) into Eq. (\ref{s06}) yields
\begin{equation}
12k^3-19k^2+8k-1 =0. \label{s08}
\end{equation}
Solving Eq. (\ref{s08}), we obtain
\begin{equation}
k \in \left \{{1\over 4}, {1\over 3}, 1\right \}. \label{s09}
\end{equation}
Equation (\ref{s07}) may be written as
\begin{equation}
{d\over dz} \left (  r+2{a'\over a} \right )  +{a'\over a} \left (
r+2{a'\over a} \right ) = k  \left (  r+2{a'\over a} \right )^2 .
\label{s010}
\end{equation}
Whence the natural change of variable
\begin{equation}
w =r +2{a'\over a}\, \cdot \label{s011}
\end{equation}
In the new variable, Eq. (\ref{s010}) becomes Bernoulli's equation
\begin{equation}
w'+{a'\over a}w =kw^2. \label{s012}
\end{equation}
The classical technique for integrating Bernoulli's equation leads
to
\begin{equation}
w  ={1\over C\, a -ka\int_0^z a^{-1}dt} \; , \label{s013}
\end{equation}
 where $C$ is an arbitrary constant of integration. Hence
 \begin{equation}
 r = {1\over C\, a -ka\int_0^z a^{-1}dt}-{2a'\over a}, \quad  k \in \left \{{1\over 4}, {1\over 3}, 1\right \}.
 \label{s014}
 \end{equation}
 To summarize, we have established the following result.

 \begin{theo}
 The operator $\mathfrak{L}$, where
 \begin{eqnarray}
 A & =&{5\,a'^2-10\,aa''\over 2a^2}\;,     \label{s015}\\
 B &=& {81\,a'^4-204\,aa'^2a''+72\,a^2a'a^{(3)}+ 84\,a^2a''^2-24\, a^2a^{(4)} \over 16 a^4},  \label{s016}
 \end{eqnarray}
is iso-spectral with the operator $\hat \mathfrak{L}$, where $\hat
A$ and $\hat B$ are defined by Eq. (\ref{f6})-(\ref{f7}), $r$ is
given by Eq. (\ref{s014}) and $ s= (A+r^2-r')/2$.
\end{theo}
\subsection{Case of the standard unit-coefficient beam ($A=0$ and $B=0$)}
According to the symmetry analysis performed in section
\ref{symanal}, this is the most symmetric instance of the
principal equation. It turns out that, in this particular case,
the principal equation can be invertibly mapped to a Chazy
equation\cite{chaz1,chaz2,chaz3}.  Indeed, Chazy's equation
\begin{equation}
y_{xxx}=2yy_{xx}-3y_x^2+\alpha (6y_x-y^2)^2,\quad \alpha
=\mbox{const.},\label{s016}
\end{equation}
admits a 3D symmetry Lie algebra spanned by the operators
\cite{clol}
\begin{equation}
Y_1={\partial \over \partial x},\quad Y_2=x{\partial \over
\partial x}-y{\partial \over
\partial y},\quad Y_3 =x^2{\partial \over
\partial x}-(2xy+6){\partial \over
\partial y}\; \cdot \label{s017}
\end{equation}
Simple calculations show that $<X_1,X_2,X_3>$ and $<Y_1,Y_2,Y_3>$
are two equivalent (up to an invertible transformation)
representations of $sl(2,\mathbb{R})$. Indeed the transformation
\begin{equation}
r=y,\quad z={2\over 3}x,\label{s018}
\end{equation}
maps $<X_1,X_2,X_3>$ into $<Y_1,Y_2,Y_3>$. The same transformation
maps the principal equation into Chazy's equation
\begin{equation}
y_{xxx}=2yy_{xx}-3y_x^2+{4\over 27} (6y_x-y^2)^2.\label{s019}
\end{equation}
The integrability of Eq. (\ref{s016}) was established by Chazy
\cite{chaz3}. The justification of the integrability of Eq.
(\ref{s016}) from a Lie symmetry standpoint is due to Clarkson and
Olver \cite{clol}.

The integration of Eq.(\ref{s016}) relies on the following
theorem.
\begin{theo}
[Chazy \cite{chaz3}] Assume that $\varphi$ and $\psi$ are two
arbitrary linearly independent solutions of the hypergeometric
equation
\begin{equation}
t(1-t){d^2\chi\over dt^2}+\left ( {1\over 2}-{7\over 6}t   \right
) {d\chi\over dt}-\sigma \chi=0,\quad \sigma={1\over
144(1-9\alpha)}\;\cdot \label{s020}
\end{equation}
Then the general solution of Eq.(\ref{s016}) is given in
parametric form by
\begin{equation}
x ={\varphi(t)\over \psi(t)},\quad y={6\over \psi(t)}{d\psi\over
dx}\,\cdot \label{s021}
\end{equation}
\end{theo}
{\bf Remark.} The case $\alpha=1/9$ was treated by Clarkson and
Olver \cite{clol}. In this degenerate case, the solution of Eq.
(\ref{s016}) is still given by Eq. (\ref{s021}), where this time,
$\varphi$ and $\psi$ are two arbitrary linearly independent
solutions of the Airy equation $\ddot \chi +(c\,t/2)\chi=0$, and
$c$ is an arbitrary constant.

For Eq. (\ref{s019}), $\sigma=-1/48$ and we may choose $\varphi$
and $\psi$ as \cite{wolf}
\begin{eqnarray}
\varphi &=& k_1\; {}_2F_1\left ( {1\over 4},-{1\over 12},{2\over
3};1-t \right )-4^{-2/3}k_2(1-t)^{1/3}\; {}_2F_1\left ( {1\over
4},{7\over
12},{4\over 3};1-t \right ), \label{s022}\\
\psi &=& k_3\; {}_2F_1\left ( {1\over 4},-{1\over 12},{2\over 3};1-t
\right )-4^{-2/3}k_4(1-t)^{1/3}\; {}_2F_1\left ( {1\over 4},{7\over
12},{4\over 3};1-t \right ), \label{s023}
\end{eqnarray}
where ${}_2F_1$ is Gauss hypergeometric function, $k_1$ to $k_4$ are
arbitrary constants, and the factor $-4^{-2/3}$ is introduced for
convenience. Since the formula (\ref{s021}) is invariant under equal
scaling of $\varphi$ and $\psi$, we may assume without loss of
generality that $k_1k_4-k_2k_3=-1$, so that Eqs.
(\ref{s022})-(\ref{s023}) involve only three arbitrary constants.
Thanks to the exact formulas \cite{vidu}
\begin{eqnarray}
{}_2F_1 \left ( {1\over 4},-{1\over 12},{2\over 3}; {\tau
(\tau+4)^3\over 4 (2\tau-1)^3} \right ) &=& (1-2\tau)^{-1/4} \label{s024}\\
 {}_2F_1 \left ( {1\over 4},{7\over 12},{4\over 3}; {\tau
(\tau+4)^3\over 4 (2\tau-1)^3 }\right ) &=& {4(1-2\tau)^{3/4}\over \tau+4}, \label{s025}\\
\end{eqnarray}
it is natural to introduce the re-parametrization
\begin{equation}
 1-t = {\tau (\tau+4)^3\over 4 (2\tau-1)^3}\,\cdot \label{s026}
\end{equation}
Thus $\varphi$ and $\psi$ reduce to
\begin{eqnarray}
\varphi & = & (k_1+ k_2 \tau^{1/3})(1-2\tau)^{-1/4} , \label{s027}\\
\psi & = &  (k_3+ k_4 \tau^{1/3})(1-2\tau)^{-1/4}. \label{s028}
\end{eqnarray}
Therefore the general solution of Eq. (\ref{s019}) is given in
parametric form by the equations
\begin{equation}
x={k_1+k_2\tau^{1/3}\over k_3+ k_4\tau^{1/3}},\quad y
={3(k_3+k_4\tau^{1/3})(3k_3\tau^{2/3}+k_4(2-\tau))\over 1-2\tau}\;
\cdot \label{s029}
\end{equation}
Solving Eq. (\ref{s029}a) for $\tau$ and substituting into Eq.
(\ref{s029}b) yield
\begin{equation}
y={3\left
[3k_3(k_3x-k_1)^2(k_4x-k_2)+2k_4(k_4x-k_2)^3+k_4(k_3x-k_1)^3 \right
]\over (k_4x-k_2)\left [2(k1-k_3x)^3-(k_4x-k_2)^3\right ]}\, \cdot
\label{s030}
\end{equation}
Finally, we obtain
\begin{eqnarray}
r &=&{6\left
[3k_3(3k_3z-2k_1)^2(3k_4z-2k_2)+2k_4(3k_4z-2k_2)^3+k_4(3k_3z-2k_1)^3
\right ]\over (3k_4z-2k_2)\left [2(2k1-3k_3z)^3-(3k_4z-2k_2)^3\right
]},  \label{s031} \\
s &= & (r^2-r')/2, \label{s033}\\
 & & k_1k_4-k_2k_3=-1.  \label{s034}
\end{eqnarray}
Thus we have established the following result.
\begin{theo}
The unit beam operator $\displaystyle{{d^4\over dz^4}}$ is
iso-spectral with the operator $\hat \mathfrak{L}$ with
\begin{eqnarray}
\hat A &=& -5r',\label{s035}\\
\hat B &=&-{1\over 4} ( 6r^{(3)}+2rr''-r^4-7r'^2), \label{s036}
\end{eqnarray}
where $r$ is defined by Eq. (\ref{s031}).
\end{theo}

\section{Conclusion}
The factorization of the canonical beam operator as the product of a
second-order linear differential operator and its adjoint results in
a system of nonlinear ordinary differential equations. We uncoupled
the system of ordinary differential equations and solved the
resulting nonlinear ordinary differential equation ({\it viz.} the
principal equation) using Lie's method. The symmetry analysis
reveals that the principal equation admits either a one-dimensional
or three-dimensional symmetry  Lie algebra. When the principal
equation admits a one-dimensional Lie algebra, its order can be
reduced  by one at most. In this case we obtained a one-parameter
family of solutions (see Eq. (\ref{s014})). For the maximally
symmetric case, we proved that the principal equation can be mapped
to a Chazy equation. The latter is solved in closed form and the
general solution (see Eq. (\ref{s031})) of the principal equation is
obtained. By reversing the order of the factorization of the
canonical beam operator, numerous nontrivial iso-spectral families
are obtained.

It is opportune to mention that Nucci \cite{nucc} calculated the Lie
point symmetries of the system (\ref{f3})-(\ref{f4}) without showing
how these symmetries are used. The truth is, for systems of ordinary
differential equations,
 successive reduction of order using symmetries can be ambiguous:
there are in general several possibilities for choosing  the new
variables as invariants \cite{waf1,waf2}. We avoided this problem by
first uncoupling the system to obtain a single nonlinear ordinary
differential equation on which we applied Lie symmetry analysis.

 \end{document}